\pdfoutput=1
\documentclass[a4paper]{IEEEtran}
\usepackage{amsmath,amsfonts,amssymb,amsthm}

\newtheorem{theorem}{Theorem}[section]

\newtheorem{definition}[theorem]{Definition}
\newtheorem{corollary}[theorem]{Corollary}

\newtheorem{property}[theorem]{Property}

\usepackage{todonotes}

\usepackage{multirow}
\usepackage{cite}
\usepackage{graphicx}
\usepackage[flushleft]{threeparttable}
\usepackage{float}
\usepackage{caption}
\usepackage{subcaption}
\usepackage{bm} 
\usepackage{dsfont}
\usepackage{wasysym}
\usepackage{nameref}
\usepackage{enumerate}
\usepackage{nomencl}
\usepackage{etoolbox}
\usepackage{colortbl}
\usepackage[capitalise]{cleveref}
\usepackage[hyphens]{url}

\newcommand{\GreenCell}{\cellcolor[RGB]{0, 200, 0}}
\newcommand{\LightGreenCell}{\cellcolor[RGB]{150, 250, 150}}
\newcommand{\RedCell}{\cellcolor[RGB]{255, 120, 120}}
\newcommand{\OrangeCell}{\cellcolor[RGB]{255, 180, 20}}

\title{Flexibility Framework with Recovery Guarantees for Aggregated Energy Storage Devices}
\date{}
\author{Michael P. Evans, Simon H. Tindemans, \IEEEmembership{Member, IEEE} and David Angeli, \IEEEmembership{Fellow, IEEE}
	\thanks{M.P. Evans completed the majority of this work while with the Department of Electrical and Electronic Engineering, Imperial College London, UK. He is currently with Octopus Energy, UK (email: michael.p.evans21@gmail.com). S.H. Tindemans is with the Department of Electrical Sustainable Energy, TU Delft, Netherlands and was also affiliated as a visiting researcher with the Alan Turing Institute, UK when the research was performed (email: s.h.tindemans@tudelft.nl). D. Angeli is with the Department of Electrical and Electronic Engineering, Imperial College London, UK, and the Department of Information Engineering, University of Florence, Italy (email: d.angeli@imperial.ac.uk). Corresponding author: S.H. Tindemans.}}
\begin{document}
\bstctlcite{MyBSTcontrol}

\IEEEoverridecommandlockouts
\IEEEpubid{\parbox{\columnwidth}{\copyright 2022 IEEE. Personal use of this material is permitted. Permission from IEEE must be obtained for all other uses, in any current or future media, including reprinting/republishing this material for advertising or promotional purposes, creating new collective works, for resale or redistribution to servers or lists, or reuse of any copyrighted component of this work in other works.} \hspace{\columnsep}\makebox[\columnwidth]{ }}

\maketitle

\IEEEpubidadjcol

\begin{abstract}
This paper proposes a framework for the procurement of flexibility reserve from aggregated storage fleets. It allows for arbitrary tree structures of aggregation hierarchy, as well as easily implementable disaggregation via broadcast dispatch. By coupling discharge and recovery modes, the proposed framework enables full-cycle capacity to be procured ahead of real time, with guaranteed recovery and exact accounting for losses. The set of feasible discharging requests is exactly encoded, so that there is no reduction in the ability to meet discharging signals, and recovery capabilities are parametrised as a single virtual battery. Included in this paper is a numerical demonstration of the construction of the constituent curves of the framework and the approach is also benchmarked against relevant alternatives.
\end{abstract}
\section{Introduction}
\label{sec:intro}
\subsection{Background}
Recent years have seen an increased proliferation of renewable energy sources onto electricity networks, leading to increased difficulty in the task of balancing supply and demand. Energy storage units offer significant potential in addressing this imbalance. For example, National Grid, the Electricity System Operator, predicts that between 28 and 75 GW of electricity storage (including vehicle-to-grid) capacity will be used to balance the British network by 2050 \cite{NG_FES}. Prior literature on balancing services has covered multiple classes of storage device, such as EV batteries \cite{Saber2012,Paterakis2016,Vandael2013,GonzalezVaya2015,Zhang2017a,Xu2016,9066844}, diesel generators \cite{Mazidi2014,Wang2015} and home storage units \cite{Paterakis2016,Tsui2012,Wang2013,6588344,8876882,8839865}; and a battery model can even capture the properties of a time-shift in electrical load \cite{7377129,Oconnell2016,9027957,9026960}. The total capabilities of units that can be modelled as \textit{virtual batteries} is therefore very significant.

A key enabling technology for large numbers of devices to participate in balancing frameworks is their \textit{aggregation}. This combination of individual units into compound entities, often termed Virtual Power Plants, offers a useful abstraction framework: in amalgamating a complex fleet into a single capability representation, one can consider there to be a single service provider. This leads to ease of interaction with the service recipient, as well as a reduction in the direct communication infrastructure connecting individual device owners and the central dispatcher.

Multiple authors have considered the aggregation of energy storage units \cite{Vandael2013,GonzalezVaya2015,Xu2016,Zhang2017a,Oconnell2016,Appino}.
The predominant form of aggregation has been the summation of power and energy limits across devices, as exemplified by \cite{GonzalezVaya2015}. While this functions as intended for homogeneous fleets, we showed in \cite{Evans2018} that such a representation only serves to outer-bound the true flexibility limits of a heterogeneous device pool. The summation-aggregation approach was improved upon in \cite{Vandael2013}, where Vandael \textit{et al.} presented a receding horizon version of the summed parameter model, with the addition of necessary constraints on total charge profiles; this same tightening of the aggregate energy constraints was then applied in \cite{Xu2016,Zhang2017a}. However, the retention of a summed power limit allowed the fleet to run at full power up to the point of full charge (equivalent to energy depletion in a discharging setting). When such an aggregation approach is utilised to compose constraints on an optimisation problem, this can still lead to the allocation of infeasible requests. Appino \textit{et al.} noted this shortfall of summation based aggregation in \cite{Appino}. While our approach is to directly rectify this problem, they instead imposed conditions so that a summation \textit{was} able to exactly capture the flexibility limits of the fleet. In doing so, these conditions necessarily restrict the aggregate flexibility, however, which may result in sub-optimal allocation when a feasible optimum is placed outside of the considered region. By contrast, the control approach presented in \cite{Evans2018} was able to exactly capture aggregate discharging capabilities. There, we presented a transform which exactly encodes fleet flexibility, under the restriction to unidirectional operation and a lack of cross-charging, termed the $E$-$p$ transform. This is the predominant tool that we will use in our analysis here, applied to coupled discharge-recharge operation. An equivalent fleet flexibility measure was presented by Cruise and Zachary in \cite{Cruise2018}. 

In addition to the aforementioned benefits, aggregation (either directly or in stages) also enables an anonymisation of individual unit owners. Moreover, where device availability is stochastic, the aggregation of units leads to a statistical smoothing of the overall behaviour, i.e. a robustification of the net output. For example, chance constraints can be applied to achieve an arbitrary risk level (as in \cite{Evans2019}) only on aggregate.

Network operators are increasingly offering frameworks for balancing services, procured through ancillary service markets (e.g. \cite{PJM,CAISO,NG_balancing,aFRR}). These markets are generally accessed by aggregators, entities which combine the capabilities of portfolios of heterogeneous units into compact delivery offers; grid operators are increasingly developing software platforms to enable scalable and automated participation via this route-to-market (see \cite{equigy} for an example). Once these markets exist, they encourage the formation of complex ecosystems in which entities of various types and sizes exchange arbitrary ancillary service provisions with one another. In this regard, a desirable feature in an aggregation scheme is the ability to abstract the capabilities of a fleet in such a way that aggregates themselves can be combined into larger aggregates. This property enables \textit{hierarchical} aggregation, which improves scalability and underpins combination and re-packaging of services. 

Examples of aggregation of heterogeneous resources  exist for various load types. For example, for thermostatically controlled loads, aggregate capability representations were presented in \cite{Trovato2015,Hao2015}. Zhao \emph{et al.}~\cite{Zhao2017} present an hierarchical aggregation method based on homothets of a prototype polytope (which encodes a standard measure of flexibility), where inner (sufficient) and outer (necessary) flexibility bounds are calculated and aggregated. Uncertainty in device parameters was introduced in \cite{9446618}, using a chance-constrained formulation. Homothet-based geometric addition was also used for active-reactive power ($P$-$Q$) flexibility in \cite{Kundu}. M\"uller \emph{et al.} introduced zonotopes as an alternative encoding of resource flexibility polytopes, based on an inner approximation (feasibility) \cite{8063901}. Zonotopes also feature straightforward aggregation, but are based on a standard set of generators instead of (homothets of) a prototype polytope. Finally, Ulbig and Andersson proposed a hierarchical aggregation scheme that features ramp rate constraints in addition to power and energy constraints \cite{ULBIG2015155}, but feasibility is no longer guaranteed after aggregation. All these methods for aggregating heterogeneous flexibility resources sacrifice either flexibility or feasibility guarantees for tractability. A qualitatively different approach is taken in \cite{9474464}, where a real-time aggregate flexibility metric is learned within a reinforcement learning framework, which is suitable for real-time control applications, but less so for markets and scheduling, as the product is not easily quantified ahead of time. Tsaousoglou \emph{et al.}~\cite{9447957}, considering the problem of optimal EV charging,  similarly forego an explicit representation of flexibility and focus on the distributed dispatch problem.

A particular challenge when deploying demand response services is the need to schedule recovery (or recharging in the case of batteries), potentially under a strict deadline. O'Connell \textit{et al.} \cite{Oconnell2016} proposed asymmetric block offers as a means of co-optimising response and rebound in a single operation. By packaging response and rebound (or equivalently discharging and recharging) into a single offer, they guarantee full recovery at a specified time, thus reducing risk exposure.  

\subsection{Contributions}
\begin{figure}
    \centering
    \includegraphics[width=\columnwidth]{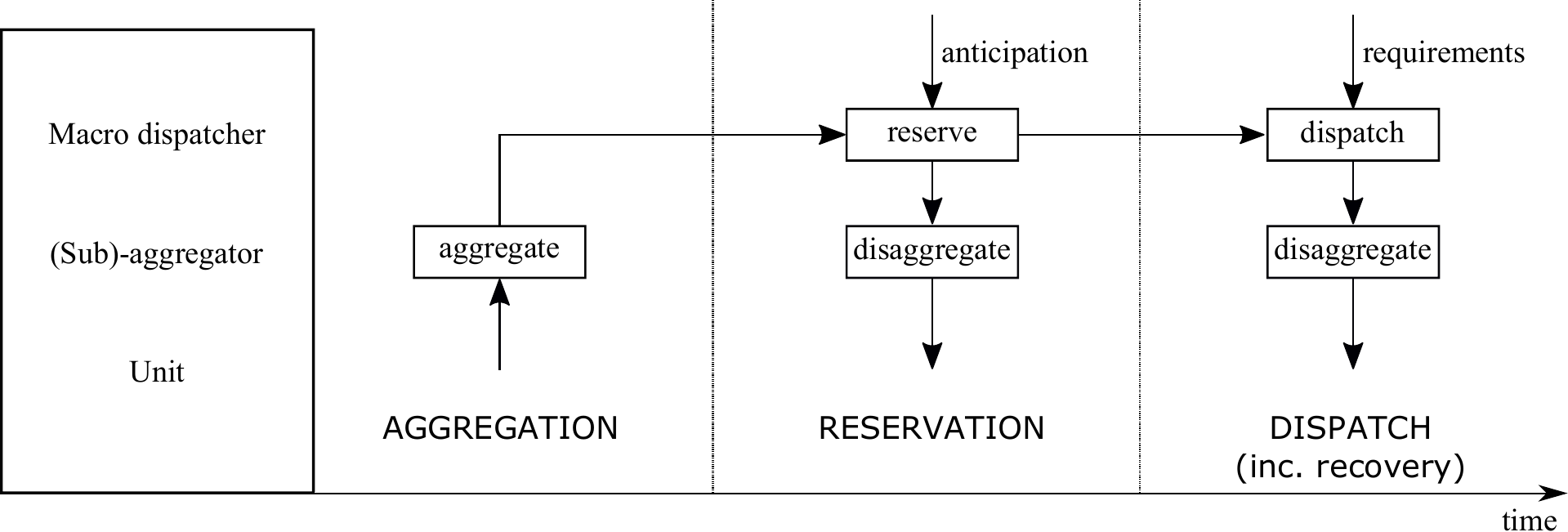}
    \caption{The separation in time of aggregation, reserve scheduling and dispatch according to the presented framework.}
    \label{fig:abstract_timing_diagram}
\end{figure}

This paper introduces a framework for hierarchical aggregation of heterogeneous storage fleets that guarantees energy recovery: the \emph{discharge-loss-recovery} (DLR) framework. This addresses the same challenge as \cite{Oconnell2016}, but extends the approach in two important ways: first, in place of a single virtual battery, the DLR representation precisely encodes the aggregate capabilities of heterogeneous storage units. Second, instead of a discrete set of offers, a continuous set of capabilities is communicated. This framework enables ahead-of-time scheduling without sacrificing achievability of recovery. As such, it enables a macro dispatcher (e.g. a system operator) to reserve flexibility in a consistent way, bounding the operational risk (and reducing associated costs) as compared to transmitting emergency requests in real time. This results in the three-stage process shown in Figure~\ref{fig:abstract_timing_diagram}. At present, there are three major use-cases for ahead-of-time reservation of flexibility: 1) grid balancing, contracted by the TSO; 2) portfolio balancing, contracted by balance responsible parties; and 3) congestion management, contracted by DSOs.

The DLR framework encodes the flexibility of storage units for the (sub-)aggregators and the macro-dispatcher, and incorporates dispatching of these units. It offers the following properties:
\begin{itemize}
\item Accounting for heterogeneity of storage units and round-trip losses.
\item Exact encoding of aggregate discharge flexibility.
\item After ensuring maximum discharge flexibility, \emph{a priori} determination of recovery requirements (energy and minimum recharge time). 
\item Hierarchical aggregation of flexibility and disaggregation via broadcast dispatch.
\end{itemize}
Reflecting the predominant use of storage under supply-shortfall conditions today, we focus on discharging requests followed by recharging recovery. The presented technique can then be applied to the opposite case via a straightforward inversion of the coordinate system used. 

The contributions of this paper are as follows:
\begin{itemize}
    \item We present a framework for coupling response and recovery capabilities that accounts for losses and asymmetry in discharging and recharging rates. This is critical to flexibility reservation involving energy-limited resources and as such improves upon the unidirectional reservation-plus-dispatch that is used by system operators for balancing today. Not only would implementing our framework ensure successful recovery, but it would also allow a system operator to control said recovery and regulate the system in doing so.
    \item Using system balancing data, we investigate the qualitative benefits that application of our framework would bring over the purchase of multiple block offers of fixed power across the delivery window.
    \item We provide closed-form equations for the aggregation of capability packets from heterogeneous device sub-fleets and describe the accompanying broadcast dispatch; together these make the approach scalable in the sense described in \cite{Vandael2013}.
    \item We qualitatively benchmark our approach against the state-of-the-art.
    \item We discuss how the framework can be applied to model demand response.
    \item We discuss how the framework can be applied in stochastic settings.
    \item We demonstrate construction of a range of DLR packets and, in doing so, show how insights can be gained from these capability encodings.
\end{itemize}

\section{Problem description}
\label{sec:prob form}
\subsection{Mathematical description}
We denote the constituent heterogeneous units of the fleet as $D_i,\ i=1,...,n$, each with extractable energy $e_i(t)$ subject to the assumed physical constraint $0\leq e_i(t)\leq \overline{e}_i$. As our focus is on the ability of a fleet to counter supply-shortfall conditions, we frame the problem based on discharging operation, choosing the sign convention that discharging power outputs are positive. We choose as our control input $u_i(t)$ the power extracted from each device, measured externally so as to take into account any inefficiencies present during discharging operation. This then leads to the following integrator dynamics:
\begin{equation}
\dot{e}_i(t)=\left\{
\begin{array}{@{}ll@{}}
-u_i(t), & \text{if}\ u_i(t)\geq0 \\
-\eta_i u_i(t), & \text{otherwise},
\end{array}\right.
\end{equation}
in which $\eta_i\in(0,1]$ denotes the round-trip efficiency of device $D_i$.

We set our problem within the paradigm of response-recovery operation, i.e. where a fleet responds by discharging available units in a nominal delivery phase, followed by a recharging phase in which the initial energy state is recovered. Initially, we focus on the discharging mode only, defining the concepts we will use for this setting specifically, before extending the same concepts to recharging mode through necessary changes in coordinates in Section~\ref{sec: symmetry to recharge}. To this end, we initially set the power limits of each device as $u_i(t)\in[0,\overline{p}_i]$. We assume a lack of cross-charging, by which we mean the use of available headroom to redistribute energy among devices.\footnote{It is shown in \cite{zachary2020scheduling} that cross-charging does not increase the feasible set for fully charged storage devices serving a unimodal signal. Moreover, cross-charging of devices with less than 100\% efficiency necessarily results in a loss of energy, which is often undesirable. } We choose as our state variable the \textit{time-to-go}, defined for each device as $x_i(t)\doteq e_i(t)/\overline{p}_i$. 
We then stack state, input and maximum power values across devices as follows:
\begin{alignat}{3}
x(t)&\doteq \begin{bmatrix}
x_1(t)\ \dots\ x_n(t)
\end{bmatrix}^T,\\
u(t) &\doteq \begin{bmatrix}
u_1(t)\ \dots\ u_n(t)
\end{bmatrix}^T,\\
\overline{p}&\doteq\begin{bmatrix}
\overline{p}_1\ \dots\	\overline{p}_n
\end{bmatrix}^T,
\end{alignat}
allowing us to rewrite the discharging dynamics in matrix form as $\dot{x}(t)=-P^{-1}u(t)$, in which $P\doteq\text{diag}(\overline{p})$. We also form the product set of the power constraints, $\mathcal{U}_{\overline{p}}\doteq[0,\overline{p}_1]\times[0,\overline{p}_2]\times...\times[0,\overline{p}_n]$,	so that our discharging constraints can be compactly written as $u(t)\in\mathcal{U}_{\overline{p}}$ and $x(t)\geq0$.

\subsection{Feasibility and flexibility}
The intention of this paper is to present a framework through which a macro dispatcher can procure flexibility reserve, ahead of real-time dispatch of the reserved capabilities. This generalises the concept of a two-stage scheduling-plus-dispatch routine (see \cite{Evans2018} for further details), since the reserve circumscribes the set of admissible schedules. For consistency, then, it is imperative that the considered framework encodes those schedules that are \textit{feasible}, by which we mean that they do not violate any physical device constraints. Initially, we focus our attention on pure discharging requests, which we allocate as positive by convention and denote\footnote{In the interest of clarity, some of the notation used here will differ from that of our prior work \cite{Evans2017,Evans2018,Evans2018a,Evans2019}} by $P^d\colon[0,+\infty)\mapsto[0,+\infty)$. We then characterise the set of feasible requests and define flexibility as follows:
\begin{definition}
	\label{def:feasible}
	The \textbf{feasible set} of discharging requests, for a system with maximum power vector $\overline{p}$ and initial state $x=x(0)$, is defined as
	\begin{equation*}
	\begin{aligned}
	\mathcal{S}_{\overline{p},x} \doteq \big\{&P^d(\cdot) \colon \exists u(\cdot),\ z(\cdot)\colon \forall t\geq 0,\ 1^Tu(t)=P^d(t),\\
	&u(t)\in \mathcal{U}_{\overline{p}},\ \dot{z}(t)=-P^{-1}u(t),\ z(0)=x,\ z(t) \geq 0 \big\}.
	\end{aligned}
	\end{equation*}
\end{definition}
\noindent By contrast, we define the \textit{flexibility} as follows:
\begin{definition}
	\label{def:flexibility}
	The \textbf{flexibility} of an energy resource is an encoding of the feasible set that has the following properties
	\begin{enumerate}
	    \item A numerical representation of admissible requests
	    \item An accompanying dispatch routine to realise admissible requests
	    \item Guaranteed feasibility of each admissible request within \emph{a priori} bounds
	\end{enumerate}
\end{definition}
\noindent We say that flexibility is maximised if the bounds on admissible signals exactly encode the feasible set, thus guaranteeing feasibility without restricting the operational envelope.

\subsection{The $E$-$p$ transform}
 We presented an example of flexibility maximisation in \cite{Evans2018}, precipitated through the restriction to pure discharging operation, termed the $E$-$p$ transform. We will utilise this here and so make the following reproductions:

\begin{definition}
	Given a discharging power request $P^d(\cdot)$, we define its E-p transform as the following function:
	\begin{equation*}
	E_{P^d}(p)\doteq\int_0^\infty \textnormal{max}\big\{P^d(t)-p,0\big\}dt,
	\end{equation*}
	interpretable as the energy required above any given power rating, $p$. The $E$-$p$ transform is convex and monotone.
\end{definition}

\begin{definition}
	\label{definition: capacity}
	We define the discharging \textbf{capacity} of a system to be the $E$-$p$ transform of the worst-case request that it can meet, i.e. \begin{equation*}
	\Omega_{\overline{p},x}(p)\doteq E_{R_{\overline{p},x}}(p),
	\end{equation*}
	with $R_{\overline{p},x}(t)\doteq\sum_{i=1}^n\overline{p}_i[H(t)-H(t-x_i)]$, in which $H(\cdot)$ denotes the Heaviside step function.
\end{definition}

\begin{property}
	\label{prop:E_p}
	A discharging request $P^d(\cdot)$ is feasible if, and only if, its $E$-$p$ transform is dominated by the capacity of the system, i.e. $P^d(\cdot)\in\mathcal{S}_{\overline{p},x} \iff E_{P^d}(p)\leq \Omega_{\overline{p},x}(p)\ \ \forall p$.
\end{property}

\begin{figure}
	\centering
	\includegraphics[width=\columnwidth]{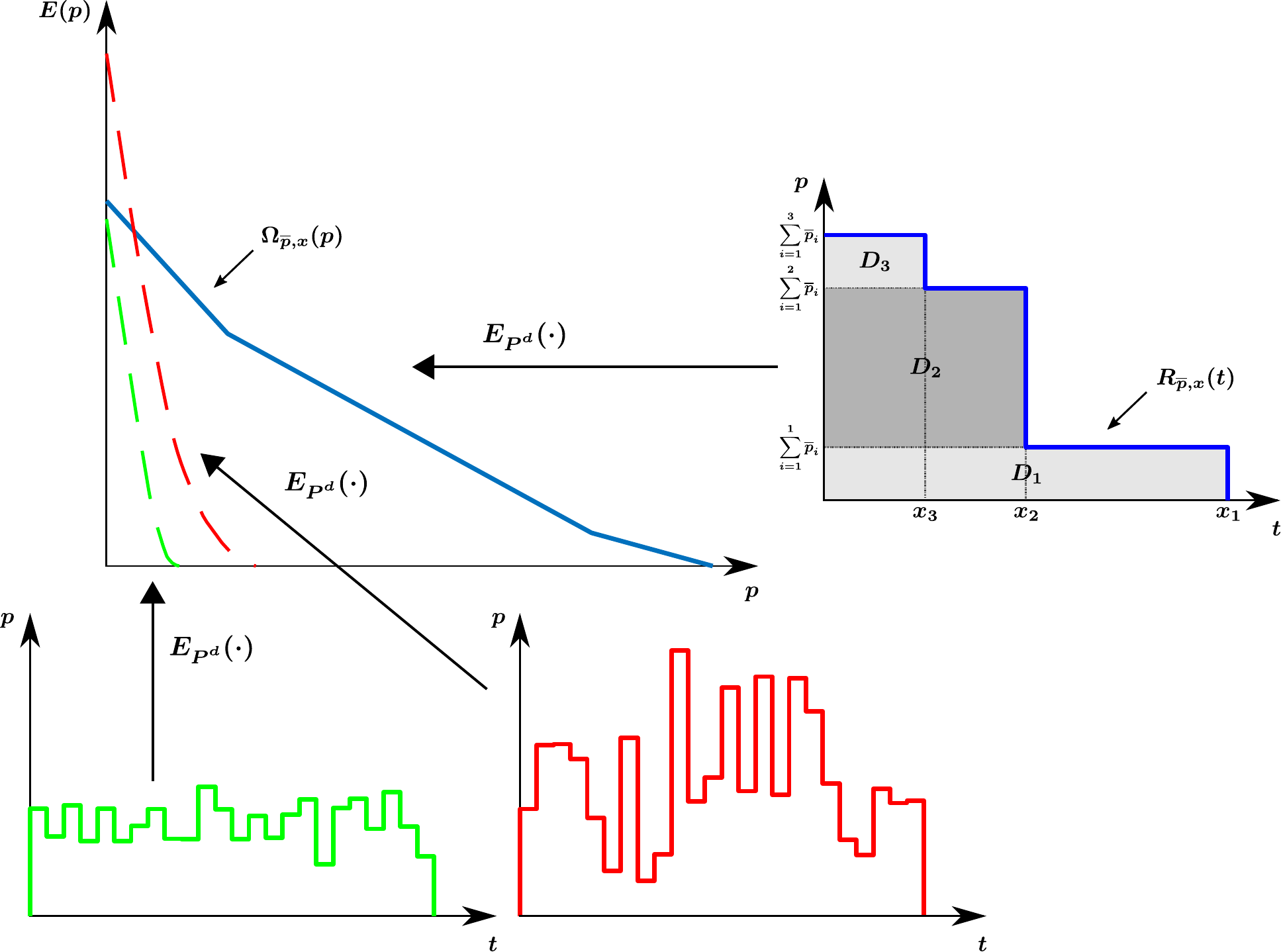}
	\caption{An example use of the $E$-$p$ transform as a binary feasibility check, implementing Property~\ref{prop:E_p}. The capacity curve is formed based on fleet parameters, and received requests are transformed into $E$-$p$ space; at which point the binary feasibility check can be made by eye. The green curve is found to be feasible, while the red curve is found to be infeasible.}
    \label{fig:Ep_example}
\end{figure}

\noindent A graphical example of the use of the $E$-$p$ transform, via Property~\ref{prop:E_p}, as a binary feasibility check can be seen in Figure~\ref{fig:Ep_example}. 
We reiterate that, under enforced conditions, the capacity curve exactly encodes the feasible set. We can therefore formally define the corresponding flexibility as follows:
\begin{definition}
	The \textbf{pure discharging flexibility} of a fleet with maximum power vector $\overline{p}$ in state $x$ is defined as
	\begin{equation*}
	\begin{aligned}
	\mathcal{F}_{\overline{p},x}&\doteq\big\{\big(p,E\big)\colon p\geq0,\ E\leq \Omega_{\overline{p},x}(p)\big\}\\
	&= \big\{\big(p,E_{P^d}(p)\big)\colon p\geq0, P^d(\cdot)\in\mathcal{S}_{\overline{p},x}\big\}.
	\end{aligned}
	\end{equation*}
\end{definition}

\noindent It also has the following property that will be useful in our analysis:
\begin{property}
	The complement to the discharging fleet flexibility, 
	$\mathcal{F}_{\overline{p},x}^\complement=\big([0,+\infty)\times[0,+\infty)\big)\setminus\mathcal{F}_{\overline{p},x}$,
	is a convex set.
\end{property}

\subsection{Default dispatch policy}
In \cite{Evans2017}, we presented a policy and showed it to maximise the flexibility of the fleet at all future time instants. We then presented a discrete time equivalent to this in \cite{Evans2018a}, which allocates a piecewise constant control input to match a piecewise constant request. In our framework, we will utilise the latter for default dispatch and therefore reproduce it here as follows. We restrict ourselves to a constant dispatch period $\Delta t$, envisioning that this will be most relevant in practical settings. Given a request that takes a constant value $P^d$ across the time interval, a target state value $\hat{z}$ is found as 
\begin{subequations}
	\begin{equation}
	\hat{z}: \sum_{i=1}^n\overline{p}_i\max\big\{0, \min\big\{x_i-\hat{z},\Delta t\big\}\big\}= P^d\Delta t,
	\label{eq:zhat_calc}
	\end{equation}
	then the input is constructed according to 
	\begin{equation}
	u_i=\overline{p}_i\cdot\textnormal{max}\bigg\{0, \textnormal{min}\bigg\{\frac{x_i-\hat{z}}{\Delta t},1\bigg\}\bigg\},
	\label{eq: discrete input}
	\end{equation}
	resulting in time-to-go updates
	\begin{equation}
	x_i \gets x_i - \frac{u_i \Delta t}{\overline{p}_i}.
	\end{equation}
	\label{eq:DTOP}
\end{subequations}
\hspace{-2mm} Note that the algorithm is applied at each time step, so explicit time-dependence is omitted. Also, unlike in \cite{Evans2018a}, we need only consider feasible requests, which leads to equality in the energy balance of \eqref{eq:zhat_calc}. 

\subsection{Symmetry to recharging operation}
\label{sec: symmetry to recharge}
Up to this point, we have considered pure discharge operation. We will now extend this to the recharging of devices between shortfall events. We denote a recharging request as $P^r\colon[0,+\infty)\mapsto(-\infty, 0]$. Analogous to time-to-go, we utilise the \textit{time-to-charge}, defined as follows:
\begin{equation}
y_i(t)\doteq\frac{\overline{e}_i-e_i(t)}{\eta_i\underline{p}_i},
\end{equation}
in which $\underline{p}_i$ denotes the (magnitude of the) maximum charge rate of device $D_i$. Equivalently, this can be obtained from a time-to-go value via the following change of (instantaneous) variables:
\begin{equation}
y_i=\frac{\overline{p}_i(\overline{x}_i-x_i)}{\eta_i\underline{p}_i},
\label{eq: x to y}
\end{equation}
in which $\overline{x}_i\doteq\overline{e}_i/\overline{p}_i$ is the maximum time-to-go value of device $D_i$. When considering pure recovery operation, then, we are able to perform the conversions $x_i\to y_i$, $P^d \to -P^r$, $\overline{p}_i\to \underline{p}_i$ and $u_i\to -u_i^r$, where $u^r$ denotes the (negative) control input while recharging, and the equivalent results hold. 

Just as the scheduling from an arbitrary state can be considered in either direction, the corresponding default dispatch policy will be given by either \eqref{eq:DTOP} or its equivalent in terms of time-to-charge. The user is alternatively able to compose the recharging policy based directly on discharging parameters as follows. Given a request that takes a constant value $P^r$ across the time interval $\Delta t$, a target (recharging) state value $\hat{y}$ is found as 
\begin{subequations}
	\begin{equation}
	\hat{y}: \sum_{i=1}^n\underline{p}_i\max\bigg\{0, \min\bigg\{\frac{\overline{p}_i(\overline{x}_i-x_i)}{\eta_i\underline{p}_i}-\hat{y},\Delta t\bigg\}\bigg\}= -P^r\Delta t,
	\end{equation}
	then the input is constructed according to 
	\begin{equation}
	u_i^r=-\underline{p}_i\cdot\textnormal{max}\bigg\{0, \textnormal{min}\bigg\{\frac{\overline{p}_i(\overline{x}_i-x_i)-\eta_i\underline{p}_i\hat{y}}{\eta_i\underline{p}_i\Delta t},1\bigg\}\bigg\},
	\end{equation}
	resulting in time-to-go updates
	\begin{equation}
	x_i \gets x_i - \frac{\eta_i u_i \Delta t}{\overline{p}_i}.
	\end{equation}
\end{subequations}

\section{The DLR framework}

\subsection{Sequential versus integrated scheduling}
Our intention is to construct a single flexibility representation that encodes the feasible set (under no cross-charging) for the full discharge-recovery cycle. This then allows for scheduling of the complete cycle ahead of real-time. 

We restrict our analysis to requests that are positive or zero (discharging) over some period, followed by negative or zero (recharging) over another period, such that full recharging of all devices is achieved; this cycle can then be repeated. We do not explicitly limit the duration of either period, but such a constraint can be communicated as an additional parameter. In line with the previous section, we shall denote a discharge-recovery request as $P^d\APLlog P^r(\cdot)$, in which the $\APLlog$ operator denotes the concatenation of two signals. 

The results in Section~\ref{sec:prob form} provide the tools to implement a two-stage solution to this problem. First, the discharging capacity curve (which, to emphasise the mode, we denote as $\Omega^d(\cdot)\doteq\Omega_{\overline{p},x}(\cdot)$) is used to determine feasibility of a discharge request $P^d(\cdot)$, and the dispatch policy \eqref{eq:DTOP} used to realise it. Then, the procedure is repeated for the recharge request $P^r(\cdot)$. However, the \emph{recharge capacity curve} $\Omega^r(\cdot)$ and, importantly, the total amount of energy required to recharge while accounting for losses, depend on the terminal state-of-charge of the storage units after the discharge operation, which in turn depends on (the E-p transform of) the discharge signal $P^d(\cdot)$. This dependency hinders the formulation of a compact flexibility representation that incorporates recharging.

\subsection{The DLR flexibility representation}

In this section, we introduce a representation of the set of feasible integrated discharge-recovery schedules that has the following properties:
\begin{itemize}
\item It maximises the discharging flexibility, so that the set of considered discharging signals exactly matches the feasible set. 
\item It exactly accounts for round-trip losses in the recharging phase. This is a crucial requirement for any aggregation method used to produce constraints on an optimisation problem. 
\item It returns recharge parameters in the form of a simple virtual battery model, characterised by recharge energy and minimum time-to-charge (equivalently recharge power and energy ratings).
\end{itemize}

As mentioned in the previous section, the recharge requirements depend on the terminal state of the discharge procedure. In order to simplify the description of recharge flexibility, we introduce a one-dimensional parameterisation of the discharging flexibility as a function of the total energy requirement $E^d$, in such a way that discharging flexibility is not sacrificed and an exact encoding of the recharge requirements is possible. This representation reduces the feasible set of recharging patterns, in exchange for a compact representation that facilitates pre-scheduling of recharging.

The proposed \emph{discharge-loss-recovery time} (DLR) representation consists of three 1-parameter functions:
\begin{itemize}
	\item [D)] The \textit{discharge} capacity curve $\Omega^d(\cdot)$ (Definition~\ref{definition: capacity});
	\item [L)] The \textit{loss} function, describing the energy required to recharge, $E^r$, as a function of the reserved discharge energy $E^d$;
	\item [R)] The \textit{recovery time} function, describing the minimum time, $y^*$, required to recharge the reserved discharge energy $E^d$.
\end{itemize}
Together, the loss and recovery functions define the $E^d$-dependent parameters of a single virtual battery for recharging. For ease of aggregation and broadcast, we choose to parameterise the loss and recovery time functions by $x^*$, the minimum time required to discharge the reserved discharge energy $E^d$, in place of $E^d$ itself. This is detailed later in the section.

\begin{figure*}[htb]
	\centering
	\includegraphics[width=\textwidth]{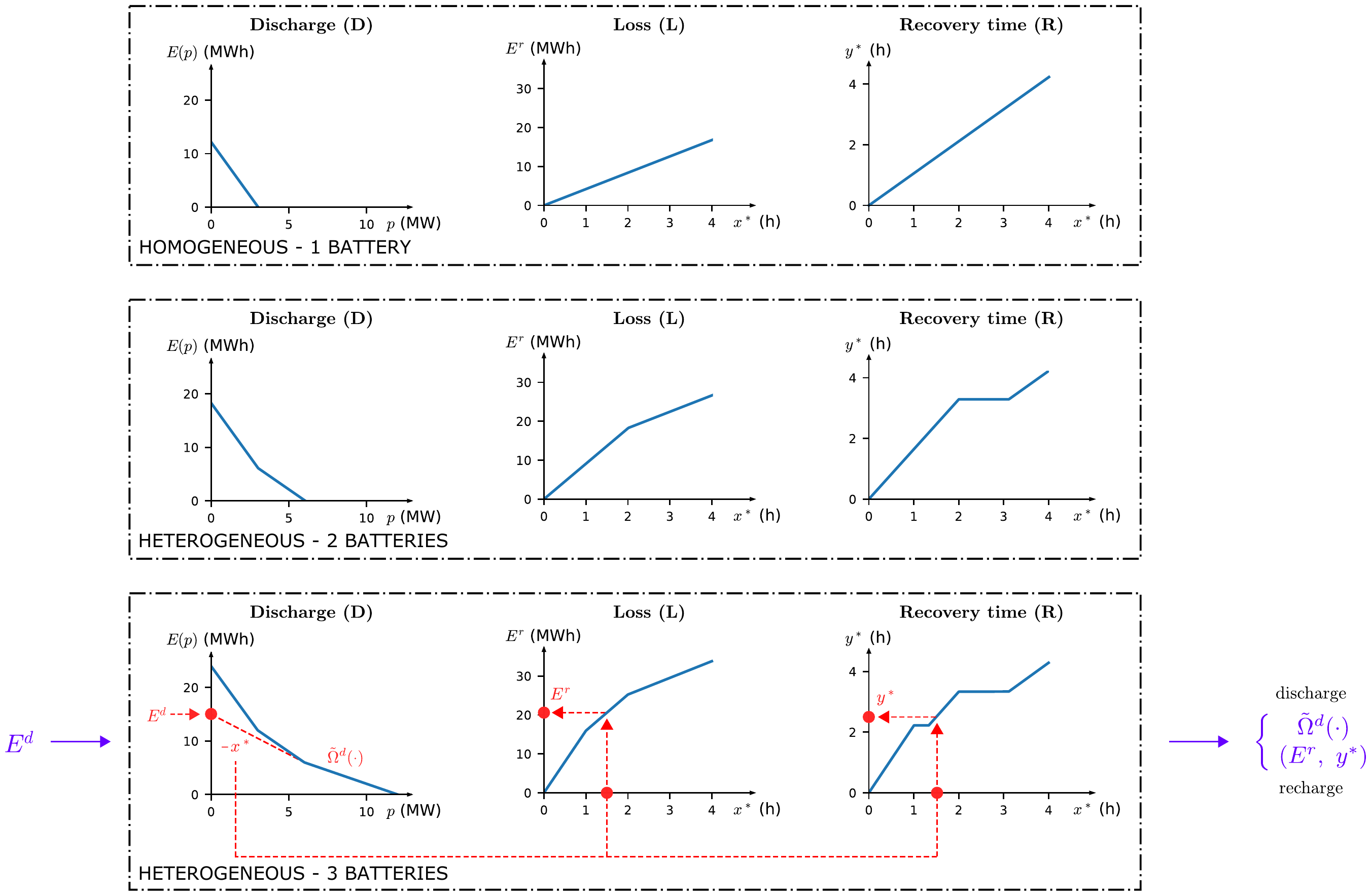}
	\caption{The DLR curves for a fleets of one (upper panel) two (centre panel) and three (lower panel) batteries. In the lower panel, red construction lines show the conversion of a discharge energy $E^d$ into the corresponding discharging and recharging flexibilities.}
	\label{fig:DLR in practice}
\end{figure*}

For a single battery, the three constituent curves of the DLR framework are composed via the following change in coordinates from the (scalar valued) parameter tuple $(\overline{p},x,\underline{p},\eta)$:
\begin{subequations}
\begin{align}
\Omega^d_{\mathrm{single}} (p) &=  x \left(  \overline{p} - p  \right), \\
E^r_{\mathrm{single}}(x^*) &= \frac{\overline{p}x^*}{\eta}, \\
y^*_{\mathrm{single}}(x^*) &= \frac{\overline{p}x^*}{\eta\underline{p}},
\end{align}
\label{eq:single_unit_DLR}
\end{subequations}
where $E^r_{single}$ and $y^*_{single}$ jointly define the recharging unit. Example curves for a single unit can be seen in Figure~\ref{fig:DLR in practice} (upper panel).

The equivalent curves for a heterogeneous fleet can be seen in Figure~\ref{fig:DLR in practice} (lower panels). Details on the construction of the DLR-curves for arbitrary fleets are provided in the remainder of this section.

\subsection{Discharging capacity and guaranteed final state}

The first element of the DLR representation, the \emph{discharging capacity} is given by the capacity curve $\Omega^d(\cdot)$, which exactly encodes the feasible set of discharging signals for a heterogeneous fleet. 

Given a feasible discharge signal $P^d(\cdot)$ and an initial state, the dispatch strategy \eqref{eq:DTOP} determines the terminal state-of-charge of each storage unit. In order to reduce the space of possible outcomes without limiting feasibility, we construct a \emph{maximum-flexibility truncated fleet} that limits the time-to-go of each storage unit. This truncated fleet is parameterised by $E^d$ and will be exactly emptied during discharge operation and refilled during recovery, \emph{without} restricting the feasible set conditioned on the discharged energy $E^d$. For a given fleet of power rating $\overline{p}$ in state $x$, we define
\begin{subequations}
	\label{eq: max-flex truncated fleet}
	\begin{equation}
	\tilde{x}_i=\min\{x_i,\ x^*\},\ i=1,...,n,
	\end{equation}
	in which the (aggregate) minimum discharge time $x^*$ is defined according to
	\begin{equation}
	x^*(E^d): \sum_{i=1}^n\overline{p}_i\min\{x_i,\ x^*\}=E^d.
	\label{eq: x star}
	\end{equation}
\end{subequations}
This truncated fleet construction satisfies the following result, the proof of which can be found in the \nameref{sec: app}:
\begin{theorem}
	\label{the:same feasibility}
	Given a request profile $P^d(\cdot)$ and a fleet with maximum power vector $\overline{p}$ in state $x$,
	\begin{equation}
	P^d(\cdot)\in\mathcal{S}_{\overline{p},x}\iff P^d(\cdot)\in\mathcal{S}_{\overline{p},\tilde{x}},
	\end{equation}
	in which $\tilde{x}$ is defined according to \eqref{eq: max-flex truncated fleet} with $E^d\geq\int_0^\infty P^d(t)dt$.
\end{theorem}
\noindent Thus, in the knowledge that $P^d(\cdot)$ will have a total energy no greater than $E^d$, the original fleet can be replaced by a truncated fleet, using \eqref{eq: max-flex truncated fleet}, and maintain its ability to meet considered signals. Moreover, if the total energy provided \emph{exactly} equals $E^d$, then the virtual fleet will be emptied and the true fleet is in a known state, regardless of the distribution of power levels in $P^d(\cdot)$.  This property gives the sought reduction in recharge complexity. 

\begin{figure}
	\centering
	\includegraphics[width=\columnwidth]{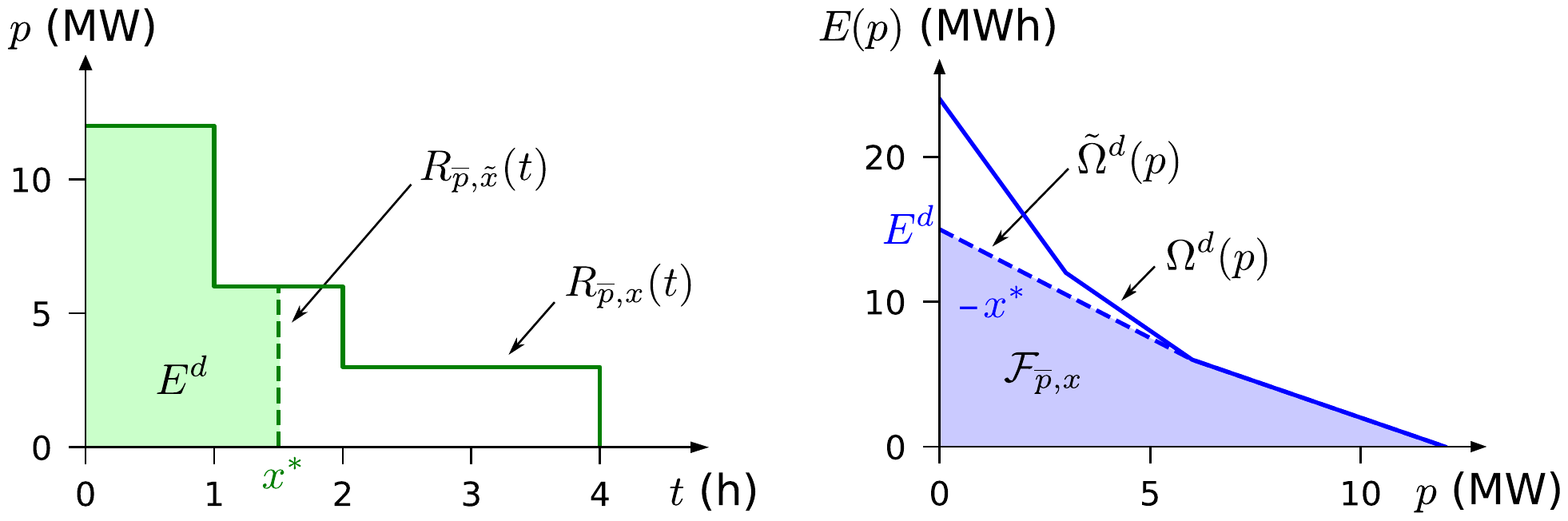}
	\caption{The correspondence between the worst-case request (left panel) and capacity curve (right panel) for a truncated fleet formed as in \eqref{eq: max-flex truncated fleet}. The green shaded area is equal to the total requested energy, while the discharging flexibility is shaded in blue.}
	\label{fig:virtual_omega}
\end{figure}

Theorem~\ref{the:same feasibility} also leads to the following Corollary:
\begin{corollary}
	\label{cor:virtual cap}
	The flexibility of a truncated fleet as in \eqref{eq: max-flex truncated fleet} can be found as
	\begin{equation}
	\mathcal{F}_{\overline{p},\tilde{x}}^\complement=\textnormal{Conv}\big(\mathcal{F}^\complement_{\overline{p},x}\cup\big\{(0,E^d)\big\}\big),
	\label{eq:conv_hull}
	\end{equation}
	in which Conv$(\cdot)$ denotes the convex hull operator.
\end{corollary}
\noindent This result gives us a means to construct the capacity curve of the truncated fleet directly from a full capacity curve, without considering constituent devices. Figure~\ref{fig:virtual_omega} demonstrates this graphically where, for simplicity of notation, we define $\tilde{\Omega}^d(\cdot)\doteq \Omega_{\overline{p},\tilde{x}}(\cdot)$. As can be seen in the figure, the flexibility of a truncated fleet can found via two equivalent methods: 1) time-truncation of the worst-case reference; 2) the convex hull operation of \eqref{eq:conv_hull} applied directly to the capacity curve. We note that, by construction, the discharge capacity curve of the truncated fleet has the property
\begin{equation}
\left. \frac{\mathrm{d} \tilde{\Omega}^d(p)}{\mathrm{d} p} \right|_{p=0^+}=-x^*. \label{eq:OmegaDerivative}
\end{equation}

\subsection{Loss and recovery time functions}

Recall that a discharge request with energy $E^d$ will exactly empty an $E^d$-truncated virtual fleet. The total energy required to recharge is therefore given by
\begin{equation}
E^r(x^*)=\sum_{i=1}^n\frac{\overline{p}_i}{\eta_i}\min\{x_i,x^*\},
\label{eq: recharge energy}
\end{equation}
with $x^*$ as in \eqref{eq: x star}. We term this the \textit{loss function} and note that it is piecewise linear, with $n$ sections in general, and concave.\footnote{The loss and recovery time functions could alternatively be parameterised by $E^d$ in place of $x^*$. The derivative of the loss function , i.e. $\frac{\mathrm{d}E^r}{\mathrm{d}E^d}$, would then be a $p$-weighted average of $1/\eta$ of participating batteries. However, we opt for $x^*$ here for increased convenience of aggregation and broadcast.}

The third and final component of the DLR representation is the \emph{recovery time function}. The initial recharge state corresponds to the truncated fleet being exactly depleted. The fastest recovery is determined by the largest time-to-charge value, which is therefore given by
\begin{equation}
y^*(x^*)\doteq\max_{i\in\{1,...,n\}}\frac{\overline{p}_i}{\eta_i\underline{p}_i}\min\{x_i,x^*\},
\label{eq: y*}
\end{equation}
with $x^*$ as in \eqref{eq: x star}. This function is monotone and piecewise linear, with $2n-1$ sections in general. 

Together, $E^r(x^*)$ and $y^*(x^*)$ describe a virtual storage unit. Acceptable recharging profiles $P^r(\cdot)$ must fully recharge this unit by supplying an amount of energy $E^r$ subject to a maximum power rating $E^r/y^*$. As it can be seen in Figure~\ref{fig:DLR in practice}, then, the constituent DLR curves are sufficient to encode, for a given total discharge energy, coupled discharging and recharging flexibilities. These guarantee feasibility of any allocated request and, as such, the DLR representation should be considered a consistent flexibility metric for integrated discharge-recovery operation. The DLR curves can also be embedded into optimisation routines to encode round-trip capabilities. 

We also point out here that, in the special case where all devices have the same round-trip efficiency and ratio of charging/discharging rates, implementation of the default discharging policy maximises \emph{charging} flexibility for a given total energy. Hence, when these parameter conditions apply, full-cycle flexibility is maximised in the sense that the recharge energy and recharge time are minimal for a given discharge energy $E^d$.

\section{Aggregation, reservation and dispatch}
\begin{figure*}[h]
	\centering
	\includegraphics[width=\textwidth]{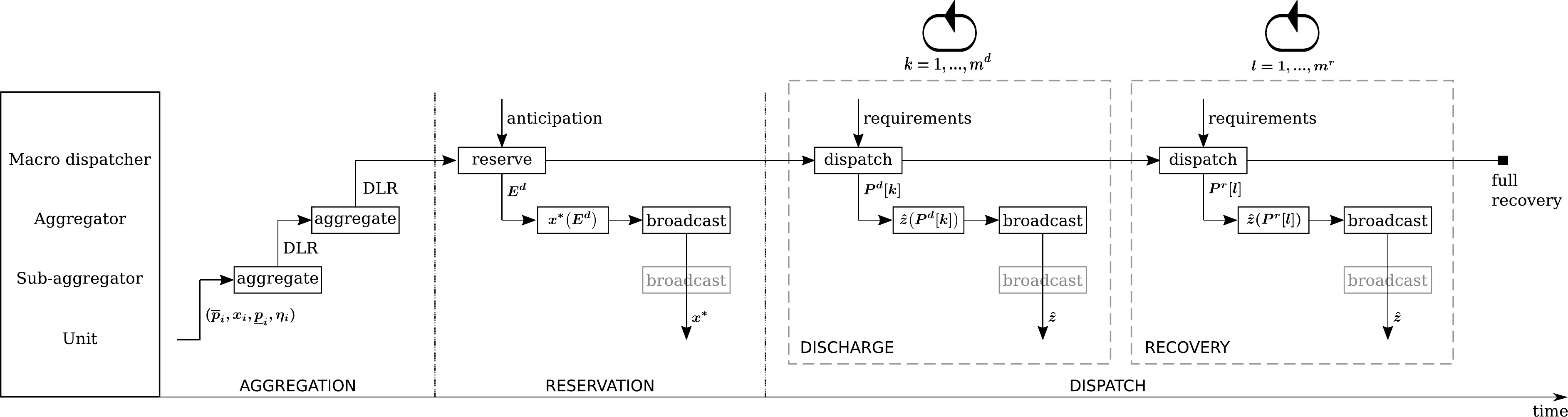}
	\caption{Aggregation, reservation and dispatch, including the messages passed between entities. $m^d\Delta t$ and $m^r\Delta t$ denote the lengths of the discharging and recharging periods, respectively. Transparent boxes show optional disaggregation steps at one or more sub-aggregator levels.}
	\label{fig:timing diagram}
\end{figure*}

The framework that we present enables the separation in time of aggregation, flexibility reservation and dispatch, as was shown in Figure~\ref{fig:abstract_timing_diagram}. Having presented the elements required for implementation of this framework, we here provide details on each of the three constituent time periods, in Figure~\ref{fig:timing diagram} and the following subsections.

\subsection{Aggregation}
The DLR representation is computed on the basis of the fleet parameter tuple $(\overline{p},x,\underline{p},\eta)$ (each vector-valued). It is then straightforward to see that DLR representations can be combined: by converting each back to fleet parameters (one virtual device per unique time-to-go value), amalgamating all of these and finally forming a single compound DLR packet. This framework therefore implicitly allows for aggregation in arbitrary hierarchies.

We would, however, like to achieve the combination of DLR packets based directly on the three constituent curves, i.e. without reverting to device parameters. Without loss of generality, we consider the aggregation of two fleets with parameters $(\overline{p}^{(1)},x^{(1)},\underline{p}^{(1)},\eta^{(1)})$ and $(\overline{p}^{(2)},x^{(2)},\underline{p}^{(2)},\eta^{(2)})$. To reduce clutter, we simplify notation by using only the subscript $1$ or $2$ to summarise the respective quantities. Each of the three DLR curves is treated in one of the following subsections.
\subsubsection{Discharge capacity curve}
We initially consider discharge flexibilities, for which we are able to derive the following result (proof in the \nameref{sec: app}):
\begin{theorem}
    \label{the:flex_addition}
    The \textbf{discharging flexibility} of an aggregate fleet can be computed by the following complementary Minkowski summation:
    \begin{equation}
        \mathcal{F}_{1+2}=(\mathcal{F}^\complement_1\oplus\mathcal{F}^\complement_2)^\complement.
    \end{equation}
\end{theorem}
\noindent This then leads directly to the following result:
\begin{corollary}
\textbf{Capacity curves} can be added as follows:
\begin{equation}
    \Omega_{1+2}(p)=\sup\{E\colon(p,E)\notin\mathcal{F}^\complement_1\oplus\mathcal{F}^\complement_2\}.
\end{equation}
\end{corollary}
The combined capacity curve $\Omega_{1+2}(p)$ also defines the fleet truncation parameter $x^*(E^d)$ through \eqref{eq:OmegaDerivative}.

\subsubsection{Loss function}
A time-to-go value $x^*$ defines all individual device recharge energy needs, regardless of their grouping into fleets. Hence the total energy requirement can be found via a simple pointwise summation, i.e.
\begin{equation}
    E^r_{1+2}(x^*)=E^r_1(x^*)+E^r_2(x^*).
\end{equation}

\subsubsection{Recovery function}
The minimum recovery time is set by the largest time-to-charge value. As such, when two curves are amalgamated, it is the maximum of these values that determines the minimum recharge time, i.e.
\begin{equation}
    y^*_{1+2}(x^*)=\max\{y^*_1(x^*),\ y^*_2(x^*)\}.
\end{equation}

Finally, we remark that the DLR representation of a fleet can be constructed in this way by aggregating DLR representations of elementary storage units, each given by \eqref{eq:single_unit_DLR}.

\subsection{Flexibility reservation}
Without loss of generality, we assume that the macro dispatcher has a direct communication link to a number of top-level aggregators. Once it has received all DLR capabilities, compound or otherwise, we envision that the macro dispatcher makes a decision as to what flexibility provision to reserve from each; and we place no restrictions on the method through which this is chosen. Once this decision has been made, the actual reservation takes place through the transmission of the amount of reserved energy (one value: $E^d$) to each top-level aggregator. Each of these then calculates the corresponding $x^*$-value, which it broadcasts among its fleet. Once reservation has been realised through the allocating of $x^*$ values, a truncated fleet is formed as each unit sets its virtual initial state $\tilde{x}_i=\min\{x_i,x^*_i\}$; this virtual device will be fully discharged and recharged.  We point out here that dispatch via the policy \eqref{eq:DTOP} ensures participation of all units, and hence so too does the allocation of devices under the proposed framework.

As an extension to this reservation method, the aggregator could in fact amend the truncation level of each device and/or the dispatch strategy, so long as the aggregate constraints were met; in this case it would broadcast a vector $x^*$ of length $n$ instead. However, this is not further investigated here.

\subsection{Broadcast dispatch}
The DLR construction is such that any admissible signal can necessarily be met by the implementation of the policy \eqref{eq:DTOP} and its equivalent for recharging. Under this combined policy, each device solely requires the $\hat{z}$ value corresponding to each request level in order to allocate its output. Moreover, this remains true regardless of the grouping of devices into fleets. As a result, the framework that we propose directly enables broadcast dispatch. There is no need to cascade signals down the tree; rather, each top-level aggregator can broadcast a single $\hat{z}$ level and the units will respond accordingly.

\subsection{Benefits of the combined framework}
We consider this approach to be \textit{scalable} (c.f. \cite{Vandael2013}), based on its combination of arbitrary aggregation hierarchy and broadcast dispatch. From their construction according to a fleet of discrete devices, each of the three DLR curves are piecewise linear. As a result, we envision that these metrics are passed between entities as vectors of their vertices; the lengths of these will scale as $\mathcal{O}(n)$.

Guaranteed recovery also means that multiple DLR packets can be sold ahead of operation for sequential use. The macro dispatcher need not fully utilise each DLR packet; in fact, it need not call upon it at all if the realised net demand does not exceed supply. Two options exist for a partial utilisation of a purchased DLR packet: 1) between discharging and recharging operation, update $E^r$ (retain $y^*$ even though it may lead to an unnecessarily conservative rate constraint); 2) before the start of the operational period, allow for scheduling updates to reduce $E^d$ and accordingly $x^*$.

\section{Numerical demonstration}
\subsection{Example construction of the DLR curves}
We now demonstrate construction of the DLR curves in a numerical example. Consider the fleet with parameters given in Table~\ref{tab:3_battery_params}, so that $n=3$ and the other parameters are such that the maximum number of linear sections is observed on each of the discharge, loss and recovery time curves: 3, 3 and 5 respectively. Firstly, consider forming the DLR curves for just Battery 1. This can be seen in Figure~\ref{fig:DLR in practice} (upper panel). Next, consider aggregating the DLR capabilities of Batteries 1 and 2. This can be seen in Figure~\ref{fig:DLR in practice} (centre panel). Finally, consider aggregating this compound DLR packet with that of Battery 3. The corresponding curves can be seen in Figure~\ref{fig:DLR in practice} (lower panel) and Figure~\ref{fig:virtual_omega}.

\begin{table}[h]
    \caption{Parameters for the three-battery example}
    \label{tab:3_battery_params}
    \centering
    \renewcommand*{\arraystretch}{1.5}
    \begin{tabular}{|c|c|c|c|c|}
      \hline
      \multirow{2}{*}{Battery \#} & \multicolumn{4}{c|}{Parameter} \\ \cline{2-5}
      & $\overline{p}_i$ (MW) & $x_i$ (h) & $\underline{p}_i$ (MW) & $\eta_i$ (\%)\\ \hline
       1&3&4&4&70 \\
       2&3&2&3&60 \\
       3&6&1&3&90 \\
       \hline
    \end{tabular}
\end{table}

\subsection{Homogeneity conditions}
\label{sec:homogeneity conditions}
We here discuss which parameter values must be uniform across the fleet for it to be considered to be homogeneous. A homogeneous fleet has discharging \textit{and} recovery curves that exactly encode the set of admissible signals and correspond to a single virtual battery.
\subsubsection{Discharging capacity curve}
The time-to-go $x_i$ dictates the slope of the capacity curve over the corresponding section. To return a curve that is a straight line therefore requires equal $x_i$ value across the fleet; this results in the discharging capabilities of a single virtual storage unit of power rating $\sum_{i=1}^n\overline{p}_i$ and total energy $\sum_{i=1}^n\overline{p}_ix_i$.

The DLR construction returns a (conditioned on $E^d$) discharging capacity curve that is a straight line between power rating $\sum_{i=1}^n\overline{p}_i$ and total energy $\sum_{i=1}^n\overline{p}_i\tilde{x}_i$.
\subsubsection{Recharging capacity curve}
Analogously to discharging mode, the requirement for a recharging capacity that corresponds to a single virtual battery is uniform time-to-charge $y_i$ across the fleet. When considering recharging between states of 0 and $x_i$ (i.e. recovery from empty to full), this condition is equivalent to uniform $(\overline{p}_i x_i)/(\eta_i\underline{p}_i)$. The homogeneous recharge capacity curve has power rating $\sum_{i=1}^n\underline{p}_i$ and total energy $\sum_{i=1}^n\overline{p}_ix_i/\eta_i$. 

When homogeneity of both discharging and recharging is enforced, the DLR construction returns a (conditioned on $E^d$) recharging capacity curve that is a straight line between power rating $\sum_{i=1}^n\underline{p}_i$ and total energy $\sum_{i=1}^n\overline{p}_i\tilde{x}_i/\eta_i$. This holds for any amount of energy dispatched.

\begin{figure}[h]
    \centering
    \includegraphics[width=\columnwidth]{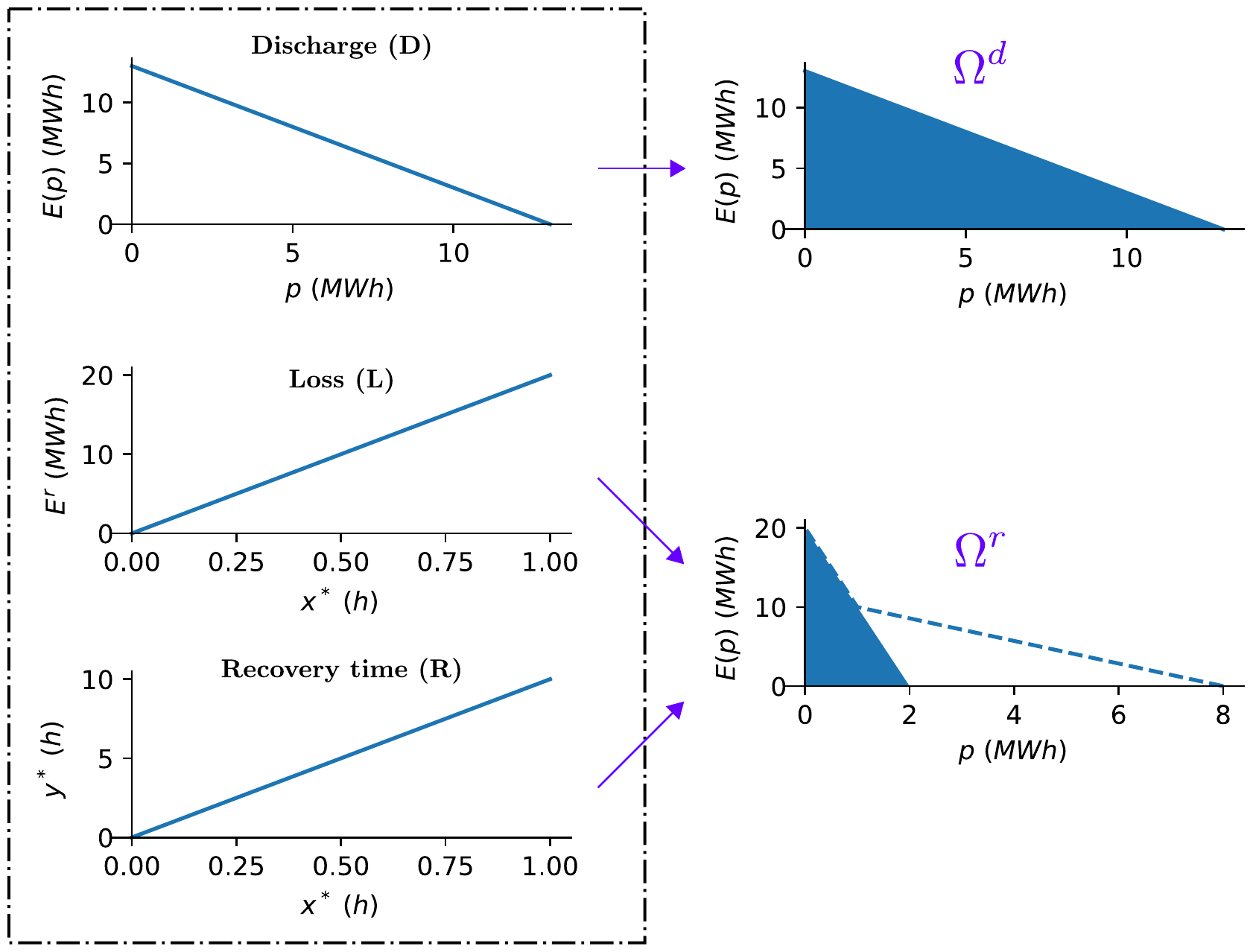}
    \caption{A two-battery example in which the fleet homogeneity conditions are met, apart from the requirement for a uniform $\underline{p}_i$ value. The DLR curves are shown on the left and used to construct discharging and recharging capacity curves on the right for $E^d=\sum_{i=1}^n\overline{p}_ix_i$. Areas below the discharging and recharging capacity curves, as returned by the DLR formulation, are shaded in blue. For comparison, the full recharging capacity curve constructed from the fleet parameters is shown in dashed blue.}
    \label{fig:2-battery example}
\end{figure}

As is to be expected, under these homogeneity conditions the loss and recovery time curves each take the form of a straight line. However, it should be noted that straight loss and recovery time curves are not sufficient conditions for homogeneity. By way of an example, Figure~\ref{fig:2-battery example} shows the DLR curves for a 2-battery fleet with parameters given in Table~\ref{tab:2_battery_params}. This is not a homogeneous fleet, hence the recharging capacity curve formed via the DLR process is a truncation of the full-fleet curve. Note, however, the straight loss curve (due to uniform $\overline{p}_i x_i/\eta_i$) and recovery time curve (due to the same condition in combination with the $\max(\cdot)$ operator).

We use this example as an opportunity to point out the priority order that has been adopted in the definition of the DLR framework. We optimise for discharging flexibility, which takes precedence over recharging efficiency. We then prioritise simplicity in encoded recharging flexibility (fully described by 2 parameters) over recharging flexibility maximisation, while guaranteeing minimum recharge time.
\begin{table}[h]
    \caption{Parameters for the two-battery example}
    \label{tab:2_battery_params}
    \centering
    \renewcommand*{\arraystretch}{1.5}
    \begin{tabular}{|c|c|c|c|c|}
      \hline
      \multirow{2}{*}{Battery \#} & \multicolumn{4}{c|}{Parameter} \\ \cline{2-5}
      & $\overline{p}_i$ (MW) & $x_i$ (h) & $\underline{p}_i$ (MW) & $\eta_i$ (\%)\\ \hline
       1&7&1&7&70 \\
       2&6&1&1&60 \\
       \hline
    \end{tabular}
\end{table}

\section{Benchmarking the DLR framework}
We here qualitatively benchmark the DLR framework against the state-of-the art in optimal scheduling-plus-operation of aggregated storage fleets. The candidate approaches are described in the following subsections, and comparisons made in Table~\ref{tab:approach_comparison}. Note that, within each broad approach, we only consider options with feasibility guarantees for practical usability. In addition to the benefits shown in Table~\ref{tab:approach_comparison}, the DLR representation is a powerful tool for simplification. This will be shown in future work due to space constraints. 

\begin{table*}[hbt!]
\caption{Benchmarking the DLR framework against alternative dispatch approaches. Green cells are preferred over orange, which are in turn preferred over red. Within the green cells, dark green denote an improvement over light green.}
\label{tab:approach_comparison}
\begin{threeparttable}
\resizebox{\textwidth}{!}{
    \centering
    \begin{tabular}{|l|c|c|c|c|c|c|c|}
    \hline
    \multirow{2}{*}{Approach} 
        & \multirow{2}{*}{discharging flexibility} & advance & \multirow{2}{*}{energy recovery} & \multirow{2}{*}{recharging flexibility} & aggregation of & visual representation & \multirow{2}{*}{dispatch method}  \\
        &  &  commitment & &  & products &  of flexibility &
         \\
    \hline
    DLR framework & \GreenCell full\tnote{*} & \OrangeCell discharge energy & \GreenCell exact, in advance & \LightGreenCell minimum time & \GreenCell yes & \GreenCell yes & \GreenCell broadcast \\ 
    Sequential $E$-$p$ dispatch & \GreenCell full\tnote{*} & \GreenCell none & \OrangeCell exact, after discharge & \GreenCell full\tnote{*} & \GreenCell yes & \GreenCell yes & \GreenCell broadcast \\
    Single virtual battery\tnote{**} & \OrangeCell limited, continuous & \GreenCell none & \RedCell not guaranteed & \OrangeCell limited, continuous & \RedCell no & \GreenCell yes & implementation dependent \\
    Asymmetric block offers & \RedCell limited, options given & \RedCell full & \GreenCell exact, in advance & \RedCell limited, options given & \RedCell no & \GreenCell yes & implementation dependent \\
    Full information scheduling & \GreenCell full & \RedCell full & \GreenCell exact, in advance & \GreenCell full & \GreenCell yes & \RedCell no & \OrangeCell direct, ahead of time \\
    Full information direct control & \GreenCell full & \GreenCell none & \GreenCell exact, full control & \GreenCell full & \GreenCell yes & \RedCell no & \RedCell direct, real time \\
    \hline
    \end{tabular}}
    \smallskip
\tiny
\begin{tablenotes}
\item[*] with the exclusion of cross-charging as discussed in \cref{{sec:prob form}}.
\item[**] parameters chosen to ensure feasibility in both directions.
\end{tablenotes}
\end{threeparttable}
\end{table*}

\subsubsection{DLR framework}
In the scheduling phase, encode the DLR curves into an optimisation problem to determine the reserved total energy $E^d$. During the operation phase, choose the aggregate power allocation in real time, subject to the equality constraint on total energy and $E$-$p$ inequality constraints (scalar during recovery), and dispatch the response power over each sample period using the DLR broadcast procedure.
\subsubsection{Sequential E-p dispatch}
During initial scheduling, reserve an $E$-$p$ curve for the fleet (optionally: parameterise this via $E^d$ and choose the maximum-flexibility truncation). Then, during discharge operation, allocate the power response over each sample period, subject to the $E$-$p$ inequality constraints only, via broadcast dispatch. After discharging operation, transform into recharging coordinates and repeat the process.
\subsubsection{Single virtual battery \cite{Vandael2013,GonzalezVaya2015,Zhang2017a,Xu2016}\protect\footnote{We consider an inner bound on parameters such that the single virtual battery encoding is guaranteed to result in the allocation of feasible signals}}
Choose a fixed set of scalar battery parameters to approximate the fleet, then incorporate this into the scheduling optimisation problem. During operation, allocate a power level for each sample period subject to the scalar constraints on power and energy.

We point out here that it is straightforward to create a feasibility-guaranteeing virtual battery from the DLR curves for a given choice of $E^d$.
\subsubsection{Asymmetric block offers \cite{Oconnell2016}}
Compute a finite set of discharge-recharge offers, each approximating the fleet as a different single virtual battery. Encode all of these offers into an optimisation problem covering the full operational period. During operation, allocate a power level for each sample period subject to the scalar constraints on power and energy.
\subsubsection{Full information scheduling \cite{Paterakis2016,Saber2012}}
Ahead of time, use full information for each device to determine an allocation of the fleet as the solution to an optimisation problem covering the full operational period. Dispatch the solution ahead of time.
\subsubsection{Full information direct control \cite{6588344,8876882}}
Omit the scheduling phase. In real time, communicate directly with each device to achieve the optimal response.

\section{Application considerations}

\subsection{Demand response}
While the DLR framework is motivated by the capabilities of a fleet of storage devices, it can be naturally applied to demand response by considering each flexible load as a \textit{virtual storage unit}. In particular, this is then well suited to the aggregation of portfolios of heterogeneous flexible loads. An additional benefit here is that the framework can be used to model a variety of load characteristics (demand turn-down, demand shifting with specific recovery requirements), via the chosen value of $\eta_i$, as can be seen in Table~\ref{tab:demand response}; at which point the heterogeneous combination can be aggregated as described. Note that, in order to model demand destruction, we extend the domain of $\eta_i$ beyond the limits of our physical battery model $\eta_i{\in}(0,1]$. The framework includes the encoding of minimum discharge ($x^*$) and recharge ($y^*$) times, but in order to ensure timely recovery of shiftable demand, a maximum allowable completion time may need to be communicated in addition.
\begin{table}[h]
    \caption{Modelling a range of demand response characteristics via the parameter $\eta_i$.}
    \label{tab:demand response}
    \centering
    \renewcommand*{\arraystretch}{1.5}
    \begin{tabular}{|l|c|}
      \hline
      Load characteristic & $\eta_i$ value \\
      \hline
      Pure demand shifting & $\eta_i=1$ \\
      Demand shifting with excess recovery & $0<\eta_i<1$ \\
      Demand shifting with incomplete recovery & $1<\eta_i<\infty$ \\
      Pure demand turn-down (no recovery) & $\eta_i=\infty$ \\
       \hline
    \end{tabular}
\end{table}

\subsection{System balancing}
The DLR framework is well suited to system balancing, in particular to markets where there are two distinct phases: reservation followed by dispatch. An example of such a market is the aFRR (automatic frequency restoration reserve) product that is common across European TSOs under the ENTSO-E umbrella organisation \cite{aFRR}. 

We focus our attention on the implementation of aFRR by TenneT TSO (Netherlands) \cite{TenneT}. It utilises a mixture of pre-purchased contracts (to ensure the ability to meet minimum requirements) and just-in-time market clearing to \textit{reserve} response, ahead of each 15~min operational period. Bids are submitted via the BRP, up to 25~min before delivery, and are separated into upwards and downwards regulation. The TSO creates a merit order list by gate closure and submits this to the aFRR platform. Instructions are then issued to \textit{dispatch} reserved flexibility within each 15~min operational period, with a response time of 30~s. The aFRR platform sends dispatch instructions to TSO-local aFRR controllers, every 1-10~s. This service is not strictly limited to large assets; in fact a number of European TSOs are developing a platform that aims to enable aggregation of small-scale resources to deliver aFRR services (Equigy, launched as a joint venture between TenneT, SwissGrid and Terna \cite{equigy}). It is conceivable that such a service may in the future also deliver local congestion management.

\begin{figure}[h]
    \centering
    \includegraphics[width=\columnwidth]{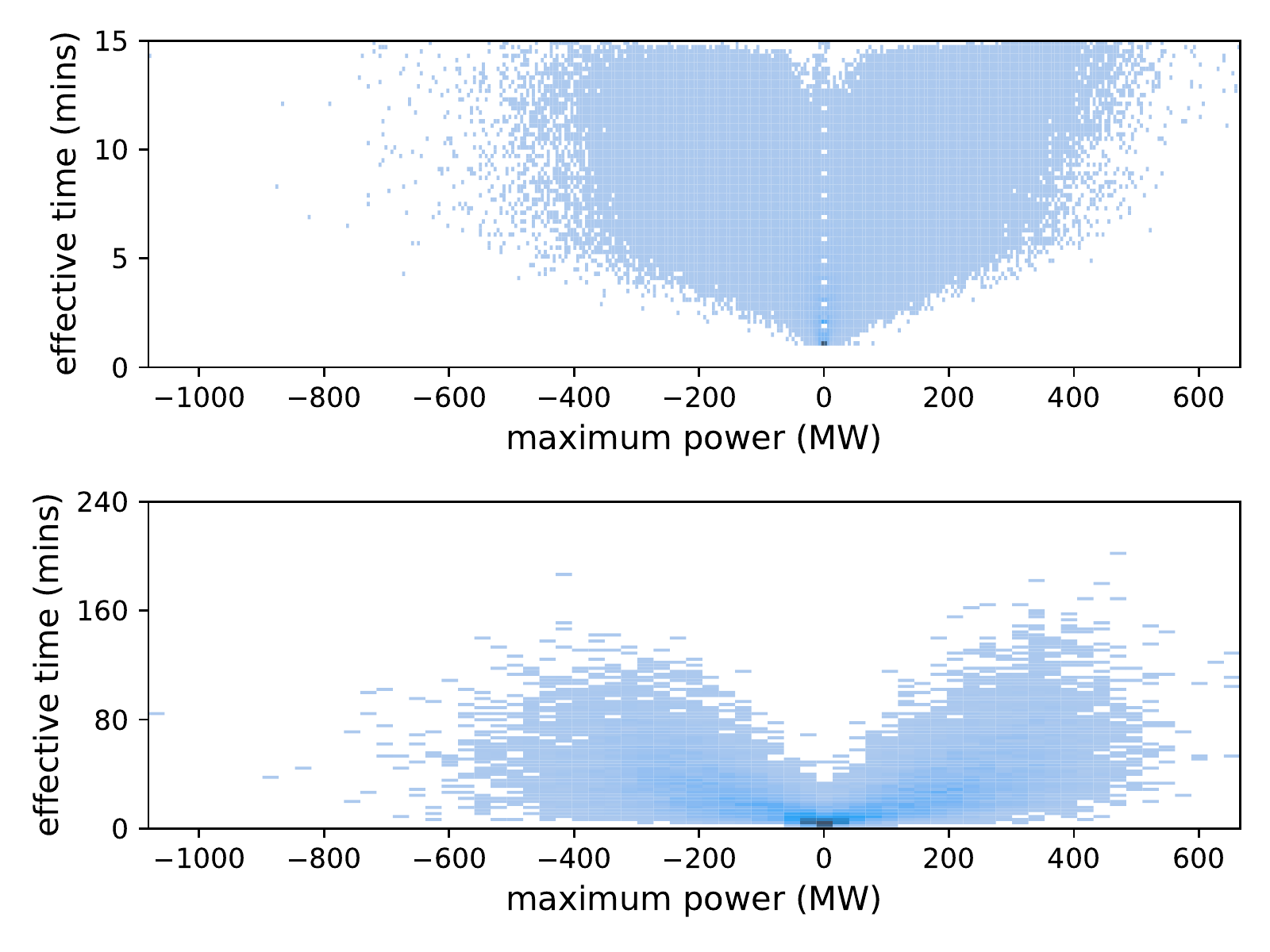}
    \caption{The distribution of energy (normalised to time) as a function of maximum power, across the Tennet dispatch dataset, for request windows of 15~min (upper panel) and 4~h (lower panel). Up-regulation requests are shown as positive and down-regulation requests as negative.}
    \label{fig:energy-power_comp}
\end{figure}

\begin{figure}[h]
    \centering
    \includegraphics[width=\columnwidth]{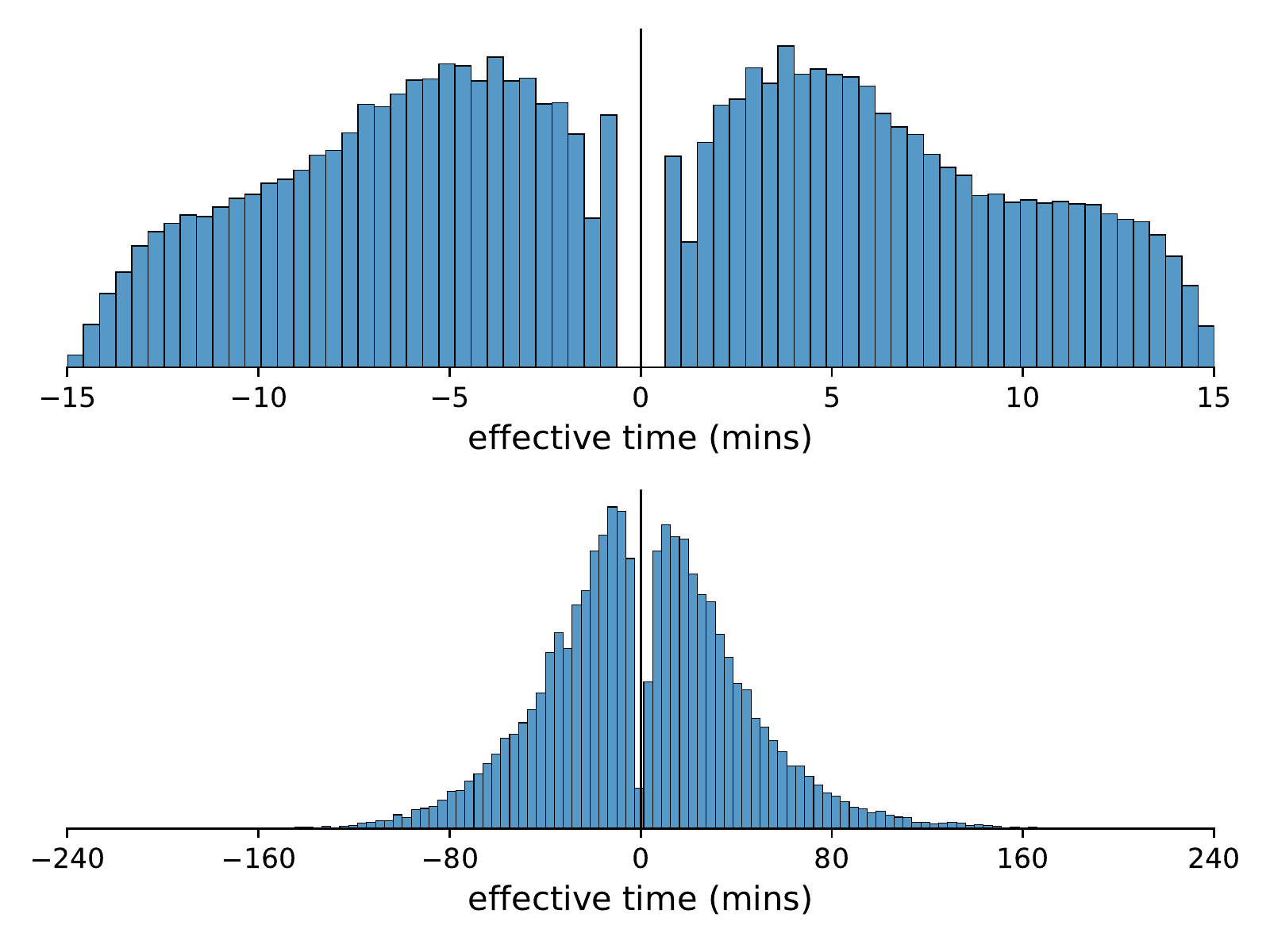}
    \caption{The distribution of effective request times, for request windows of 15~min (upper panel) and 4~h (lower panel). Times corresponding to up-regulation requests are shown as positive and down-regulation as negative.}
    \label{fig:E-p_distribution}
\end{figure}

In order to assess the suitability of the DLR framework for application to the TenneT Netherlands aFRR market, we analyse the reported dispatch data for the 5~y period covering 2017 through 2021, at 1~min granularity, publicly available from the TenneT website \cite{TenneT_data}. Currently, the market requires block bids of constant power per PTU (15 min). We will investigate the possible benefits of more flexible products for such a service. To do so, we split the 5 y traces into request segments, at 15 min as well as 4 h granularity to assess longer-duration requirements (e.g. for pre-contracted services). We analyse up- and down-regulation separately, because these are treated as independent products; for each we neglect periods that have the value 0 MW throughout their duration. Firstly, we investigate the distribution of dispatched energy versus the maximum power level within each time window. The results of this can be seen in Figure~\ref{fig:energy-power_comp}, where the vertical axis shows the dispatched energy divided by maximum power -- a duration that is necessarily capped at the window size. Figure~\ref{fig:E-p_distribution} shows the projection of these densities on the effective-time axis, separately for up- and down-regulation.

As can be seen in these figures, it is very rare for the observed peak power to be required across the entire settlement period; rather, a large range of effective time spans is observed. This is true across all peak power levels and the effective time tends to increase with maximum power magnitude (correlation coefficients are shown in Table~\ref{tab:corr_coeffs}). While for 15~min observation windows the effective time can take values up to the full duration (i.e. a constant power requirement), the 4~h requests do not exceed an effective time of 203~min and the majority have an energy requirement corresponding to less than half of the period length.

\begin{table}[h]
    \caption{Correlation coefficients between maximum power and effective time for the results of Figure~\ref{fig:energy-power_comp}. Negative values result from the sign convention that down-regulation requests are negative.}
    \label{tab:corr_coeffs}
    \centering
    \renewcommand*{\arraystretch}{1.5}
    \begin{tabular}{|c|c|c|}
      \hline
      \multirow{2}{*}{Window size} & \multicolumn{2}{c|}{Service direction} \\ \cline{2-3}
      & Up-regulation & Down-regulation\\ \hline
       15~min &0.64&-0.59 \\
       4~h &0.62&-0.59 \\
       \hline
    \end{tabular}
\end{table}

\begin{figure}[h]
    \centering
    \includegraphics[width=\columnwidth]{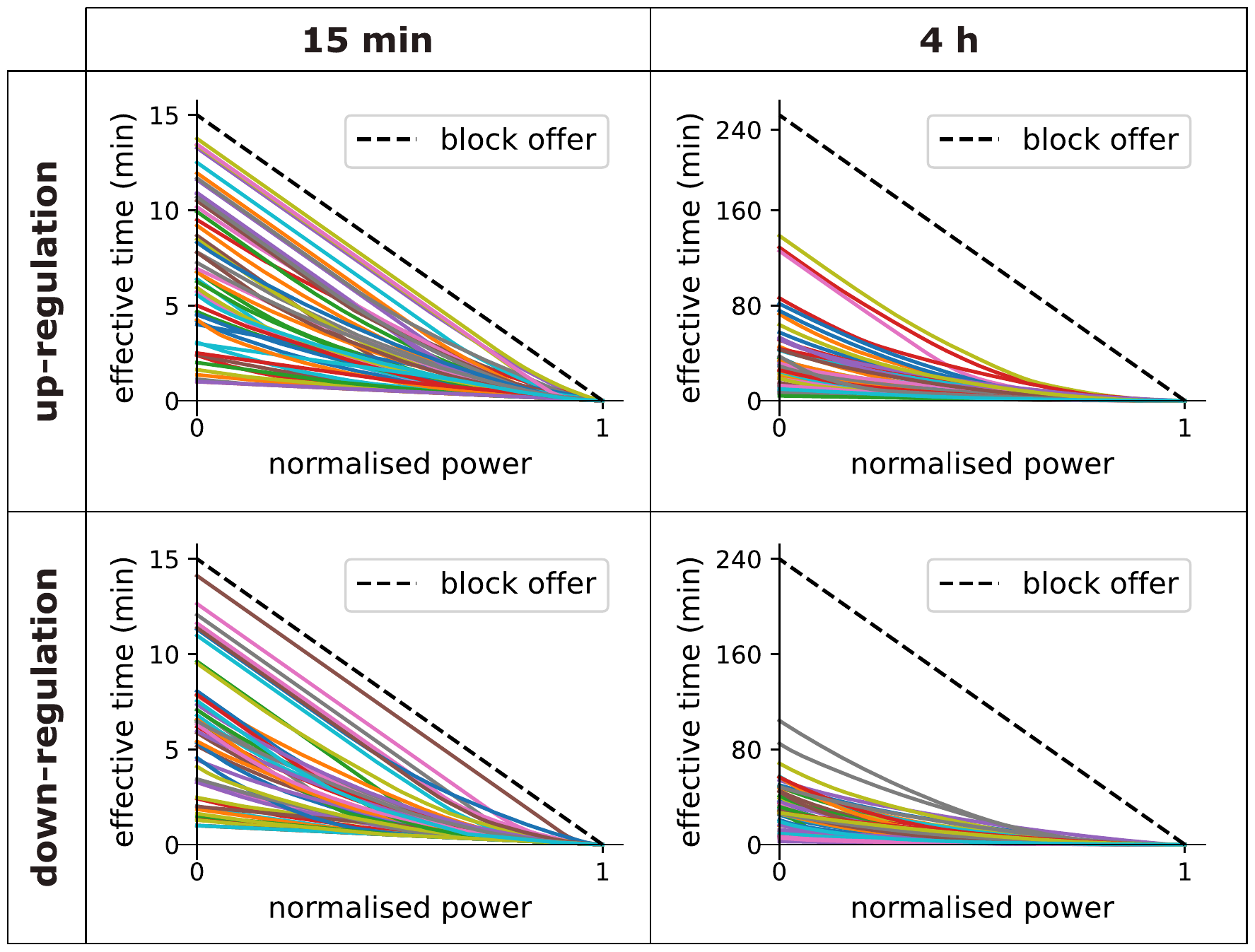}
    \caption{Discharge capacity curves composed from the Tennet dispatch request dataset. 50 random samples are shown for each of up-regulation at 15~min granularity (upper left panel), down-regulation at 15~min granularity (lower left panel), up-regulation at 4~h granularity (upper right panel) and down-regulation at 4~h granularity (lower right panel).}
    \label{fig:E-p_samples}
\end{figure}

These results suggests that, under the currently used framework, the TSO is required to purchase a wide variety of different power-energy products, given the restriction to block offers in which the energy is given by 15~min multiplied by the fixed power. We propose an improved alternative to this: that it instead purchases a single (compound) DLR product that is a better match for its needs for the period. To demonstrate the potential benefits of this approach over the restriction to block offers of fixed duration, we plot the discharge capacity curve across 50 randomly-sampled dispatch signals from the TenneT dataset in Figure~\ref{fig:E-p_samples}. We do this in both the up- and down-regulation directions, and for each sample normalise the request so that its maximum power is equal to 1; in this way we can visually compare across different request magnitudes.

As it can be seen in Figure~\ref{fig:E-p_samples}, the discharge capacity curves lie below the straight line corresponding to a single block offer, covering more of this space in the 15~min case than in the 4~h case. By construction, the distribution in Figure~\ref{fig:E-p_distribution} corresponds to the distribution of the y-axis intercepts of these (rescaled) $E$-$p$ transforms. The combination of these results shows that, for 15 min request windows, the majority of requests use less than the full energy that is reserved in the block offer. For 4 hour request windows, the mismatch between power and energy is exacerbated. This suggests that advance balancing contracts could make allowances for energy constraints and, for example, be provided by duration-limited batteries. 

To summarise, use of the discharge capacity curve in place of (a stack of diverse) block offers allows the system operator to procure the exact reserve that it requires, in one well-defined and transparent interaction with aggregators. Use of the DLR approach captures this benefit and in addition extends it into settings of response followed by recovery. The utilisation of energy-limited assets (physical batteries or virtual batteries formed of, for example, flexible loads) for balancing necessitates such a recovery mode. This is challenging for system balancing because activation will either prohibit rapid recovery or cause uncontrolled recovery (which causes further balancing concerns). The DLR approach lets the balancing authority explicitly control round-trip operation and dispatch the energy recovery process. This may result in additional efficiency by reducing the need to purchase regulation products in the opposing direction. 

\subsection{Dealing with uncertainties}

As presented thus far, the DLR framework is a deterministic methodology. A DLR packet is computed to encode the fleet capabilities and compared to system needs, for example via embedding into an optimisation problem. The same framework can also be applied in stochastic settings, by taking a scenario approach and performing this comparison for each scenario. This leads to two possible application settings: 1) a range of DLR packets is computed, sampling uncertainty in the capabilities of the fleet, and a determination made on the ensemble of DLR curves (stochastic capabilities); 2) a single DLR packet is computed and compared to an ensemble of different system requirements (stochastic system needs). A combination of the two would also be possible, in which both the capabilities and system needs were sampled from distributions.

By way of an example, we demonstrated chance-constrained optimisation under stochastic capabilities for unidirectional operation, using the capacity curve alone, in \cite{Evans2019}, via a quantile approach.
The same approach could be used to make determinations on stochastic fleet capabilities for response-recovery operation using the DLR curves. As the constituent curves (previously the discharging curve alone) serve to encode the flexibility constraints in any given scenario, this approach leads to a consistent deduction of the stochastic flexibility bounds across the range of scenarios sampled.

\section{Conclusion and future work}
This paper has presented a framework for the procurement of reserve flexibility with guaranteed recovery. The DLR representation encodes full discharging capabilities and easily interpretable recovery capabilities that exactly account for losses. Feasibility is guaranteed, so that any request deemed admissible based on the purchased reserve can necessarily be met by the fleet. We have provided an example construction of the constituent curves and benchmarked our approach against the state-of-the-art.

In future work, the authors intend to extend these results into settings where cross-charging is allowed to take place. In particular, they aim to investigate the applicability of the DLR framework in such circumstances and characterise any resulting reduction in optimality. In addition, the possibility to use simplified fleet representations (which lead to even more compact DLR representations) will be investigated.

\section*{Acknowledgements}
This work was supported by EPSRC studentship 1688672, by The Alan Turing Institute under the EPSRC grant EP/N510129/1 and by the Lloyd’s Register Foundation-Alan Turing Institute programme on Data-Centric Engineering under the LRF grant G0095. The authors would also like to thank the Isaac Newton Institute for Mathematical Sciences for support and hospitality during the programme Mathematics of Energy Systems (EPSRC grant number EP/R014604/1). Special thanks to Pierre Pinson for enlightening discussions on the topic.

\section*{Appendix}
\renewcommand{\theequation}{A.\arabic{equation}}	
\setcounter{equation}{0}
\label{sec: app}
\subsection{Proof of Theorem~\ref{the:same feasibility}}
Given $\tilde{x}$, define $\tilde{p}\doteq\sum_{i\colon x_i\geq \tilde{x}}\overline{p}_i$. Using this, we can write down the capacity curve based on $\tilde{x}$ as follows:
\begin{equation}
\Omega_{\overline{p},\tilde{x}}(p)=\left\{
\begin{array}{@{}ll@{}}
\Omega_{\overline{p},x}(p),& \text{if}\ p\geq \tilde{p} \\
\Omega_{\overline{p},x}(\tilde{p})+(\tilde{p}-p)\tilde{x}, & \text{otherwise}.
\end{array}\right.
\label{eq: omega x'}
\end{equation}
In particular, $\Omega_{\overline{p},\tilde{x}}(0)=E^d$. We can therefore deduce that
\begin{equation}
\left.\begin{array}{@{}cc@{}}
E_{P^d}(p)\leq\Omega_{\overline{p},x}(p)\ \ \forall p \\
E_{P^d}(0)\leq E^d
\end{array}\right\}\iff E_{P^d}(p)\leq \Omega_{\overline{p},\tilde{x}}(p)\ \forall p,
\end{equation}
and so the result follows from Property~\ref{prop:E_p}.
\qed

\subsection{Proof of Theorem~\ref{the:flex_addition}}
Consider two fleets of devices: $(\overline{p}^{(1)},x^{(1)})$ of size $n$, indexed by descending time-to-go values (i.e. $x_i \ge x_{i+1}$) and $(\overline{p}^{(2)},x^{(2)})$ that consists of a single device. With slight abuse of notation, we denote the parameters of the latter by $\overline{p}^{(2)},x^{(2)}$. We define $j: x^{(1)}_j \ge x^{(2)} > x^{(1)}_{j+1}$, with $x^{(1)}_0\doteq x^{(2)}$ and $x^{(1)}_{n+1}\doteq 0$. Now, consider the resultant capacity curve of the  combined fleet, as follows:
	\begin{equation}
	\Omega_{1+2}(p)=\left\{
	\begin{array}{@{}ll@{}}
	\Omega_1(p)+\overline{p}^{(2)}x^{(2)}, & \text{if}\ p\leq p' \\
	\Omega_1(p-\overline{p}^{(2)}), & \text{if}\ p\geq p'+\overline{p}^{(2)}\\
	\Omega_1(p')+(p'+\overline{p}^{(2)}-p)x^{(2)}, & \text{otherwise},
	\end{array}\right.
	\label{eq: summation omega}
	\end{equation}
	in which $p'\doteq \sum_{1=1}^{j}\overline{p}^{(1)}_i$. Denote the flexibilities corresponding to the two fleets as $\mathcal{F}_1$ and $\mathcal{F}_2$. Denote their vertices as the sets $\mathcal{V}_1$ and $\mathcal{V}_2$, which are given by
	\begin{subequations}
	\begin{alignat}{2}
	\mathcal{V}_1&=\mathcal{C}(\Omega_1)\cup\{(0,+\infty),(+\infty,0)\},\\
	\mathcal{V}_2&=\{(0,+\infty),(0,\overline{p}^{(2)}x^{(2)}),(\overline{p}^{(2)},0),(+\infty,0)\},
	\end{alignat}
	\end{subequations}
	in which $\mathcal{C}(\Omega)$ denotes the vertices of a capacity curve $\Omega(\cdot)$. We can then compute the Minkowski sum as follows:
	\begin{equation}
	\begin{aligned}
	    \mathcal{F}^\complement_1\oplus\mathcal{F}^\complement_2&=\textnormal{Conv}(\{a+b\colon a\in\mathcal{V}_1,\ b\in\mathcal{V}_2\})\\
	    &=\textnormal{Conv}(\{a+b\colon a\in\mathcal{C}(\Omega_1),\ b=(0,\overline{p}^{(2)}x^{(2)}),(\overline{p}^{(2)},0)\}\\
	    &\hspace{14mm}\cup\{(0,+\infty),(+\infty,0)\}),
	    \end{aligned}
	\end{equation}
	in which Conv$(\cdot)$ denotes the formation of the convex hull. This returns the corresponding set of vertices
	\begin{equation}
	\begin{aligned}
	    &\{(\sum_{i=1}^k\overline{p}^{(1)},\sum_{i=k}^n\overline{p}^{(1)}_k x_k+\overline{p}^{(2)}x^{(2)}),\ k=0,...,j-1\}\\
	    &\quad \cup\{(\sum_{i=1}^k\overline{p}^{(1)}+\overline{p}^{(2)},\sum_{i=k}^n\overline{p}^{(1)}_k x_k),\ k=j,...,n\}\\
	    &\quad \cup\{(0,+\infty),(+\infty,0)\}\\&=\mathcal{C}(\Omega_{1+2})\cup\{(0,+\infty),(+\infty,0)\},
	    \end{aligned}
	\end{equation}
	so the result follows for the addition of a single device to a complex fleet. The general result then follows because the Minkowski summation operator is associative.
\qed

\bibliographystyle{IEEEtran}
\bibliography{bib_all,biblio_manual3}

\begin{IEEEbiographynophoto}{Michael P. Evans}
received the BA and MEng degrees in Engineering from the University of Cambridge, UK, in 2015 and the PhD in Electrical and Electronic Engineering from Imperial College London, UK, in 2019. He is currently with Octopus Energy, London, UK.

His research interests lie in the optimal control of power systems, with a particular focus on the provision of aggregate flexibility from large fleets of heterogeneous assets.
\end{IEEEbiographynophoto}

\begin{IEEEbiographynophoto}{Simon H. Tindemans (Member, IEEE)} received the M.Sc. degree in physics from the University of Amsterdam, Amsterdam, The Netherlands, in 2004 and the Ph.D. degree in theoretical biophysics from Wageningen University, Wageningen, The Netherlands, in 2009. From 2010 to 2017, he was a Research Associate, Marie-Curie Intra-European Fellow and Research Fellow with the Control and Power Group, Imperial College London, U.K. He is currently an Assistant Professor with the Department of Electrical Sustainable Energy, Delft University of Technology, Delft, The Netherlands.

His research concerns the use of statistical approaches to model, analyze and control the power system. Within this domain, he has a speciﬁc research interest in the aggregate behavior of large numbers of heterogeneous and unpredictable resources, and the ability to characterize and control their aggregate flexibility.
\end{IEEEbiographynophoto}

\begin{IEEEbiographynophoto}
{David Angeli (Fellow, IEEE)} received the B.S. degree in Computer Science Engineering and the Ph.D. in Control Theory from University of Florence, Italy, in 1996 and 2000, respectively.

Since 2000 he was an Assistant and Associate Professor (2005) with the Department of Information Engineering, University of Florence and since 2008, he joined the EEE Department of Imperial College London, where he is currently a Professor in Nonlinear Network Dynamics and the Director of the MSc in Control Systems. He has been a Fellow of the IEEE since 2015 and a Fellow of IET since 2018. He is the recipient of the 2021 International Honeywell Medal from the Institute of Measurement and Control.

He authored more than 100 journal papers in the areas of Stability of nonlinear systems, Control of constrained systems (MPC), Chemical Reaction Networks theory and Smart Grids.
\end{IEEEbiographynophoto}
\end{document}